\documentclass[makeidx]{article}
\title{
Sharp bounds for the first eigenvalue of a fourth order Steklov problem
\footnote{Classification AMS $2000$: 58J50, 35P15, 35J40 \newline 
Keywords: Fourth order Steklov problem, Eigenvalues, Harmonic functions,  Lower bounds
\newline $^1$ The second author acknowledges financial supported by GNSAGA and MURST of Italy}}
\author{Simon Raulot and Alessandro Savo$^1$}
\date{\today}
\usepackage{makeidx}
\usepackage{amssymb , amsmath, amsthm}
\newtheorem{defi}{Definition} 
\newtheorem{thm}[defi]{Theorem}

\newtheorem{rem}[defi]{Remark}
 \newtheorem{prop}[defi]{Proposition}

\newtheorem{cor}[defi]{Corollary}
 
\newcommand{\twosystem}[2]{\left\{\begin{aligned} &#1\\ &#2\end{aligned}\right.}
\newcommand{\threesystem}[3]{\left\{ \begin{aligned}&#1\\ &#2\\&#3\end{aligned}\right.}

\newcommand{\parte}[1]{\smallskip\noindent {\rm#1)}\,\,}

\newcommand{\abs}[1]{\lvert{#1}\rvert}

\newcommand{\reals}{{\bf R}}
\newcommand{\sphere}[1]{{\bf S}^{#1}}
\newcommand{\hyp}[1]{{\bf H}^{#1}}
\newcommand{\real}[1]{{\bf R}^{#1}}

\newcommand{\bd}{\partial}
\newcommand{\starred}[1]{#1^{\star}}

\newcommand{\derive}[2]{\frac{\bd #1}{\bd#2}}
\newcommand{\isoper}{\dfrac{\abs{\bd\Omega}}{\abs{\Omega}}}

\begin{document}

\maketitle

\begin{abstract} 
We study the biharmonic Steklov eigenvalue problem on a compact Riemannian manifold $\Omega$ with smooth boundary.
We give a computable, sharp lower bound of the first eigenvalue of this problem, which depends 
only on the dimension, a  lower bound of the Ricci curvature of the domain, a lower bound of the mean curvature of its boundary and the inner radius. The proof is obtained  by estimating the isoperimetric ratio of non-negative subharmonic functions on $\Omega$, which is of independent interest. We also give a comparison theorem for geodesic balls. 
\end{abstract}
\large²


\section{Introduction}\label{intro}



\subsection{The biharmonic Steklov problem} \label{biharmonic}


Let $\Omega$ be a compact, connected Riemannian  manifold of dimension $n$ with smooth boundary $\bd\Omega$. We consider the following fourth order eigenvalue boundary problem
\begin{equation}\label{biharmonic}
\twosystem
{\Delta^2 f=0\quad\text{on}\quad\Omega}
{f=0, \quad \Delta f=q\derive fN\quad\text{on}\quad \bd\Omega,}
\end{equation}

where $N$ is the inward  unit normal and $\Delta=-{\rm div}\nabla$ is the Laplacian defined by the Riemannian metric of $\Omega$; the sign convention is that, in Euclidean space, $\Delta=-\sum_j \bd^2/\bd x_j^2$. 

The first eigenvalue of \eqref{biharmonic} has the following variational characterization:
\begin{equation}\label{var}
q_1(\Omega)=\inf\Big\{\dfrac{\int_{\Omega}(\Delta f)^2}{\int_{\bd\Omega}\big(\derive fN\big)^2}: f=0\quad\text{on}\quad \bd\Omega
\Big\}.
\end{equation}
This problem was introduced by Kuttler and Sigillito \cite{KS} and Payne \cite{P}.
For a review of the main facts about $q_1(\Omega)$, also when the boundary is not smooth, we refer to \cite{BFG} :in this paper, we always assume smoothness of $\bd\Omega$. For a recent lower bound of $q_1$, see \cite{wangxia}. It is known that $q_1(\Omega)$ is positive and simple, and that any first eigenfunction $f$ does not change sign on $\Omega$ (see \cite{Ku}). Moreover, it turns out that $q_1(\Omega)$ has the following interesting characterization in terms of harmonic functions:
\begin{equation}\label{betatwo}
q_1(\Omega)=\inf\left\{\dfrac{\int_{\bd\Omega}h^2}{\int_{\Omega}h^2}: \text{$h$ is harmonic on $\Omega$}\right\}.
\end{equation}
The infimum in \eqref{betatwo} is attained precisely when $h=\Delta f$, where $f$ is a first eigenfunction of \eqref{biharmonic} (see \cite{Ku}).
In particular, taking $h=1$ one observes the following upper bound by the isoperimetric ratio:
\begin{equation}\label{isoper}
q_1(\Omega)\leq\isoper
\end{equation}
where $|\Omega|$ (resp. $|\partial\Omega|$) denotes the volume of $\Omega$ (resp. of $\partial\Omega$) for the induced Riemannian measure. It turns out that, if $\Omega$ is a geodesic ball in a simply connected manifold of constant sectional curvature, then equality holds in \eqref{isoper} (see Theorem \ref{main} below). 

\smallskip

The main scope of this paper is to give sharp lower bounds for $q_1(\Omega)$. Let us review some known results. Isoperimetric inequalities for plane domains were given in \cite{Ku}, \cite {KS} and \cite {P}. In \cite{P}, Payne considers convex domains in $\real n$ and shows that, if $w_{\Omega}$ is the minimal distance between two parallel hyperplanes which enclose $\Omega$, then
\begin{equation}\label{payne}
q_1(\Omega)\geq\dfrac{2}{w_{\Omega}}.
\end{equation}
So, convex Euclidean domains which are "thin" have large first eigenvalue. We will prove below that this is a quite general principle which applies to a large class of manifolds. 

\medskip

For Riemannian manifolds, Wang and Xia prove in \cite{wangxia} that, if $\Omega$ has nonnegative Ricci curvature, and $\bd\Omega$ has mean curvature bounded below by $H>0$ (hence positive everywhere), then
\begin{equation}\label{wx}
q_1(\Omega)\geq nH,
\end{equation}
with equality if and only if $\Omega$ is a Euclidean ball.  We should mention that inequality \eqref{wx} is also a consequence of a lower bound of $q_1(\Omega)$ by the first eigenvalue of a certain Steklov problem for differential forms, recently proved by the authors (see Theorem $10$ in \cite{RS}).

\smallskip

Wang and Xia proof of \eqref{wx} makes use of the Reilly formula: this approach, however, seems to be hard to apply when the curvature of the domain (or the mean curvature of its boundary) assume negative values somewhere. 

\smallskip

In this paper, using a Laplacian comparison argument, and a symmetrization procedure developed in \cite{S}, we give a sharp, explicit lower bound of $q_1(\Omega)$: see Theorem \ref{main}. The method produces a positive lower bound for all compact Riemannian manifolds with boundary, and in particular it  applies also in negative curvature. Moreover, our bound implies (and actually improve) both \eqref{payne} and \eqref{wx}: see Remarks \ref{remwx} and \ref{rempayne}.

\smallskip

Finally, we prove a comparison theorem for geodesic balls, analogous to  Cheng comparison theorem for the first Dirichlet eigenvalue: see Theorem \ref{cheng}.


\subsection{The main estimate}\label{me}


Let $S$ be the shape operator of $\bd\Omega$ relative to the inner unit normal $N$: recall that if $X$ is a vector tangent to $\bd\Omega$, then $S(X)=-\nabla_XN$, where $\nabla$ is the Levi-Civita connection of $\Omega$. The mean curvature $\mathcal{H}$ of $\bd\Omega$ is the function defined by
$$
\mathcal{H}=\dfrac{1}{n-1}{\rm tr}S,
$$
and the sign convention is that $\mathcal{H}=1$ for the boundary of the unit ball in $\real{n}$.

\begin{defi}\label{cb} 
Let $K$ and $H$ be arbitrary real numbers. We say that the 
$n$-th dimensional domain $\Omega$ has curvature bounds $(K,H)$ if:

\smallskip

$-$ {\it the Ricci curvature of $\Omega$ is bounded below by $(n-1)K$},

$-$ {\it the mean curvature of $\bd\Omega$ is bounded below by $H$}.
\end{defi}

Let $R$ be the {\it inner radius} of $\Omega$, defined as
the maximal radius of a ball included in $\Omega$. Clearly:
$$
R=\max\{{\rm dist}(x,\bd\Omega): x\in \Omega\}.
$$
We then prove that $q_1(\Omega)$ is uniformly bounded below by a positive constant depending only on $K,H,n,R$.  To make this constant explicit, introduce the function:
\begin{equation}\label{sk}
s_K(r)=\threesystem
{\dfrac {1}{\sqrt{K}}\sin (r\sqrt{K})\quad\text{if}\quad K>0,}
{r\quad\text{if}\quad K=0,}
{\dfrac {1}{\sqrt{\abs{K}}}\sinh (r\sqrt{\abs{K}})\quad\text{if}\quad K<0,}
\end{equation}
and let
\begin{equation}\label{theta}
\Theta (r)=(s'_K(r)-Hs_K(r))^{n-1}.
\end{equation} 
Note that $\Theta$ depends on $K$ and $H$, and that $\Theta(0)=1$.
In Proposition \ref{kasue} it will be shown that $\Theta$ is positive on $[0,R)$, and that $\Theta(R)=0$ if and only if $\Omega$ is a ball in the simply connected manifold $M_K$ of constant curvature $K$ (we will call $M_K$  a {\it space form}).

\smallskip

Here is the main estimate.
\begin{thm}\label{main} 
Assume that $\Omega$ has curvature bounds $(K,H)$ (see Definition \ref{cb}). Let $\Theta$ be the function defined  in \eqref{theta}. Then:
\begin{equation}\label{mainbound}
q_1(\Omega)\geq\dfrac{1}{\int_0^R\Theta(r)\,dr}.
\end{equation}
The inequality is sharp: it reduces to an equality when $\Omega$ is a ball in the space form $M_K$ of constant curvature $K$  (in this case, the right-hand side is the isoperimetric ratio $\abs{\bd\Omega}/\abs{\Omega}$).
\end{thm}

The proof will be given in Sections \ref{mainproofone} and \ref{mainproof}. It turns out that, besides balls in $M_K$, equality holds also for other domains, like flat cylinders (see Corollary \ref{ircor}); however, we will not study the complete equality case of Theorem \ref{main}. The proof of Theorem \ref{main} relies on an estimate of the isoperimetric ratio of non-negative subharmonic functions, which is of independent interest (see Theorem \ref{boundaryintegral}). 

\smallskip

In Section \ref{ee} we have explicited a number of lower estimates of $q_1(\Omega)$, which follow directly from \eqref{mainbound}.

\smallskip

Let us see some other consequences of Theorem \ref{main}. For simplicity we assume $K\in\{0,1,-1\}$ so that the  space form $M_K$ is, respectively, the Euclidean space $\real n$, the sphere $\sphere{n}$ and the hyperbolic space $\hyp n$. 
\begin{cor}\label{maincor} 
Assume that $\Omega$ has curvature bounds $(K,H)$, and that:

\parte a $H>0$ if $K=0$;

\parte b $H\in\reals$ if $K=1$;

\parte c $H>1$ if $K=-1$.

\smallskip

Let $\bar\Omega$ be the unique ball in $M_K$ having boundary of (constant) mean curvature $H$. Then
$$
q_1(\Omega)\geq q_1(\bar\Omega)=\dfrac{\abs{\bd\bar\Omega}}{\abs{\bar\Omega}},
$$
with equality if and only if $\Omega=\bar\Omega$.
\end{cor}

For the proof, see Section \ref{proofcor}.
\begin{rem}\label{remwx} 
\rm   Corollary \ref{maincor} can be seen as a generalization of  the Wang-Xia bound \eqref{wx}, which is the statement in (a): in fact, the Euclidean ball of mean curvature $H$ has first eigenvalue $q_1(\bar\Omega)=\abs{\bd\bar\Omega}/\abs{\bar\Omega}=nH$. However, the main bound \eqref{mainbound} is actually stronger and a straightforward calculation gives 
$$
q_1(\Omega)\geq \dfrac{nH}{1-(1-RH)^n}\geq nH
$$
because $1-RH\geq 0$, with equality only for Euclidean balls (see Theorem \ref{eeone}a). 
\end{rem}

Corollary \ref{maincor}b applies to spherical domains; in particular we get that, if $\Omega$ is a domain
in $\sphere n$ having boundary of non-negative mean curvature, then the ball $\bar\Omega$ is just the hemisphere and one has:
\begin{equation}\label{hemi}
q_1(\Omega)\geq q_1(\Bar\Omega)=\dfrac{2\abs{\sphere{n-1}}}{\abs{\sphere n}}.
\end{equation}
with equality if and only if $\Omega=\bar\Omega$. 

\smallskip

Finally,  Corollary \ref{maincor}c implies that, for all  domains in ${\bf H}^n$ with mean curvature bounded below by $1$, one has the simple bound
\begin{equation}\label{mckean}
q_1(\Omega)> n-1,
\end{equation}
and the equality is asymptotically approached by the hyperbolic ball of radius $r$, as $r\to\infty$: this is seen immediately by estimating the isoperimetric ratio of hyperbolic balls (but see also Theorem \ref{eetwo}).

\smallskip

The lower bound \eqref{mckean} is a kind of Mc Kean inequality for the eigenvalue $q_1(\Omega)$
(the original Mc Kean inequality states that the first eigenvalue of the Laplacian for the Dirichlet boundary conditions satisfies the bound $\lambda_1(\Omega)\geq (n-1)^2/4$ for {\it all} domains in $\hyp n$).


\subsection{Lower bounds by the inner radius}\label{innerradius}


Perhaps, a more interesting feature of the main lower bound \eqref{mainbound}
is its dependence on the inner radius. We first illustrate this fact on a special case
(see Theorem \ref{eeone}b for a proof).
\begin{cor}\label{ircor} 
Assume that $\Omega$ has non-negative Ricci curvature and that $\bd\Omega$ has non-negative mean curvature.  Then:
$$
q_1(\Omega)\geq\dfrac1R.
$$
Equality holds for flat cylinders, that is, for  all Riemannian products 
$\Omega=N\times [0,2R]$ where $N$ is any closed manifold and $R>0$. 
\end{cor}

\begin{rem}\label{rempayne} 
\rm Corollary \ref{ircor} applies to mean-convex (in particular, convex) Euclidean domains, and in that case it improves the bound \eqref{payne}. In fact, just observe that, if the  strip between two parallel planes contains $\Omega$, it contains a ball of maximal radius inside $\Omega$. Hence $w_{\Omega}\geq 2R $ which implies
$$
q_1(\Omega)\geq\dfrac1R\geq \dfrac{2}{w_{\Omega}}.
$$
However, very often one has $1/R> 2/w_{\Omega}$: for example, if $\Omega$ is close to the equilateral triangle circumscribed to the unit disk one has $w_{\Omega}\sim 3$
while $2R\sim 2$.
\end{rem}

We then observe the following  rough, general estimate.
\begin{cor}\label{rough} 
Let $\Omega$ be a domain with curvature bounds $(K,H)$. Assume that $R\leq 1$. Then:
$$
q_1(\Omega)\geq\dfrac{c}{R}.
$$
for a positive constant $c$ depending only on $K,H$ and $n={\rm dim}\,\Omega$.
\end{cor}

For the proof, let $\bar R\in (0,\infty]$ be the first positive zero of $\Theta$: then we know that $R\leq \bar R$ (see Proposition \ref{kasue}).  Let $C$ be the maximum value of $\Theta$ on the interval $[0,\min\{1,\bar R\}]$ (note that $C$ depends only on $K,H$ and $n$). Since
$$
\int_0^R\Theta(r)\,dr\leq RC,
$$
the corollary follows from Theorem \ref{main} by taking $c=1/C$. 

\smallskip

In conclusion, if the Ricci curvature of $\Omega$ and the mean curvature of $\bd\Omega$ are {\it uniformly bounded below},  then $q_1(\Omega)$ becomes larger and larger as the inner radius $R$ tends to zero. 
 
With that in mind, the following general principle holds:
\begin{rem}\label{thin}
Thin domains have large first eigenvalue.
\end{rem}

Finally, we remark that a lower bound of the mean curvature is essential in order to have a positive, uniform lower bound of $q_1(\Omega)$.  In fact, let $B_r$ be the ball of radius $r$ in the unit sphere $\sphere{n}$.
As $r\to \pi$ (the diameter of the sphere) we see that $\abs{\bd B_r}$ tends to zero while the volume of $B_r$ approaches the volume  of the sphere. The isoperimetric ratio $\abs{\bd B_r}/\abs{B_r}$ tends to zero, and, by \eqref{isoper}, so does $q_1(B_r)$. Note however that the mean curvature of  $\bd B_r$ tends to $-\infty$ as $r\to\pi$.


\subsection{An upper bound of Cheng type}\label{ct}


Let $M$ be a manifold with Ricci curvature of $M$ bounded below by $(n-1)K$, with $K\in \reals$, and  $B(x_0,r)$ be the geodesic ball in $M$ with center $x_0$ and radius $r$. Cheng proves in \cite{C} the following comparison theorem for the first eigenvalue of the Laplacian under Dirichlet boundary conditions:
$$
\lambda_1(B(x_0,r))\leq\lambda_1(B_K(r)),
$$
where $B_K(r)$ is any ball of radius $r$ in $M_K$, the space form of constant curvature $K$.  It turns out that the same fact holds for $q_1$. Precisely:

\begin{thm}\label{cheng} 
Assume that the Riemannian manifold $M$ has Ricci curvature bounded below by $(n-1)K$. Then, for all $x_0\in M$ and for all $r$ less than the injectivity radius of $M$ at $x_0$, we have:
$$
q_1(B(x_0,r))\leq q_1(B_K(r)),
$$
where $ B_K(r)$ is any ball of radius $r$ in $M_K$.
\end{thm}
The assumption on the injectivity radius is done to insure that the ball $B(x_0,r)$ has smooth boundary. For the proof, see Section \ref{chengproof}.


\subsection{Boundary integral estimates}\label{bie}


The proof of Theorem \ref{main} is obtained by estimating the isoperimetric ratio (that is, the quantity $\int_{\bd\Omega}h/\int_{\Omega} h$) of any non-negative subharmonic function $h$ on $\Omega$.  Precisely, the following estimate holds.

\begin{thm}\label{boundaryintegral} 
Let $\Omega$ be a compact domain with smooth boundary.  Assume that $\Omega$ has curvature bounds $(K,H)$. If $h$ is any non-trivial, non-negative subharmonic function on $\Omega$ (that is, $h\geq 0$ and $\Delta h\leq 0$ on $\Omega$), then
\begin{equation}\label{smvl}
\dfrac{\int_{\bd\Omega}h}{\int_{\Omega}h}\geq \dfrac{1}{\int_0^R\Theta(r)\,dr},
\end{equation}
where $\Theta$ is defined in \eqref{theta}. 

\smallskip

If $\Omega$ is a ball in the space form $M_K$ and $h$ is harmonic, then equality holds; in that case \eqref{smvl} reduces to the classical mean-value lemma for harmonic functions on balls:
$$
\dfrac{\int_{\bd\Omega}h}{\int_{\Omega}h}=\isoper.
$$
\end{thm}
The proof of Theorem \ref{boundaryintegral} is given in Section \ref{mainproof}.

In Theorem $3.1$ of \cite{GS}, P. Guerini and the second author give a lower bound of 
$\int_{\bd\Omega}h/\int_{\Omega}h$ when $h$ is non-negative and satisfies $\Delta h\leq \mu h$ for $\mu$  less than a suitable positive constant.
However, this bound is sharp only when $K=H=0$.

\smallskip

Let us sketch the idea of the proof of \eqref{smvl}. Introduce the function $F:[0,\infty)\to\reals$ defined as
\begin{equation}\label{fh}
F(r)=\int_{\{\rho>r\}}h,
\end{equation}
where $\{\rho>r\}$ denotes the set of points of $\Omega$ at distance greater than $r$ from the boundary.
As $h$ is nonnegative and subharmonic, the function $F$ is shown to satisfy the following simple differential inequality:
$$
F''-\dfrac{\Theta'}{\Theta}F'\geq 0
$$
in the sense of distributions. Integrating this inequality twice one gets the desired bound. All this will be explained in detail in the next section. 
This symmetrization procedure has been used in \cite{S} and \cite{GS}  to get spectral estimates in various contexts.


\subsection{Proof of Theorem \ref{main}}\label{mainproofone} 


Let $f$ be a first eigenfunction of problem \eqref{biharmonic}, and set $g=\Delta f$. Then $g$ is harmonic on $\Omega$ and on $\bd\Omega$ one has $g=q_1(\Omega)\bd f/\bd N$. Now:
$$
\int_{\Omega}g^2=\int_{\Omega}(\Delta f)^2=q_1(\Omega)\int_{\bd\Omega}\Big(\derive fN\Big)^2=\dfrac 1{q_1(\Omega)}\int_{\bd\Omega}g^2,
$$
that is
$$
q_1(\Omega)=\dfrac{\int_{\bd\Omega}g^2}{\int_{\Omega}g^2}.
$$
As $g$ is harmonic, the function $h=g^2$ is subharmonic and non-negative. Then,  applying Theorem \ref{boundaryintegral}  to $h$ we immediately obtain the assertion:
$$
q_1(\Omega)\geq \dfrac{1}{\int_0^R\Theta(r)\,dr}.
$$
Finally, if $\Omega$ is a geodesic ball in $M_K$, then the right-hand side of the previous inequality is the isoperimetric ratio $\abs{\bd\Omega}/\abs{\Omega}$ (see \eqref{density}). On the other hand, by \eqref{isoper}, one has also $q_1(\Omega)\leq 
\abs{\bd\Omega}/\abs{\Omega}$, hence equality holds.


\section{Proofs}\label{mainproof}



\subsection{Known facts on the distance function}\label{basic}


The proof of Theorem \ref{boundaryintegral} is based on the calculation of the second derivative, in the distributional sense, of the function $F$ defined in \eqref{fh}. This calculation was done in \cite{S} and involves the Laplacian of the distance function to the boundary: see Proposition \ref{calc} below. In this section  we state the results we need from \cite{S}; for convenience of the reader, we recall the main arguments of the proof  in the Appendix below. 

\smallskip

Let $\Omega$ be a compact domain with smooth boundary  and let
$\rho:\Omega\to [0,\infty)$ be the distance function to the boundary:
$$
\rho(x)={\rm dist}(x,\bd\Omega).
$$
Then $\rho$ is only Lipschitz regular, and $\abs{\nabla\rho}=1$ almost everywhere on $\Omega$.  The set where $\rho$ is singular (i.e. not $C^1$) is called the {\it cut-locus} and is denoted by ${\rm Cut}(\bd\Omega)$: it is closed, and has measure zero in $\Omega$.  Since the boundary of $\Omega$ is assumed smooth, the function $\rho$ is smooth on a strip near $\bd\Omega$, precisely on $\rho^{-1}\big[0,{\rm inj}(\bd\Omega)\big)$ where ${\rm inj}(\bd\Omega)={\rm dist}(\bd\Omega,{\rm Cut}(\bd\Omega))$ is the {\it injectivity radius} of the normal exponential map. 

\smallskip

We define the distributional Laplacian of $\rho$ in the usual way:
$$
(\Delta\rho,f)\doteq \int_{\Omega}\rho\Delta f,
$$
for all smooth functions $f$ compactly supported in  the interior of $\Omega$, where $(\cdot,\cdot)$ denotes the duality between a test-function and a distribution.  
 We want to estimate $\Delta\rho$ from below. To that end, assume that $\Omega$ has Ricci curvature bounded below by $(n-1)K$,  that $\bd\Omega$ has mean curvature bounded below by $H$ and let $\Theta$ be the function defined in \eqref{theta}. 
Then one has the following bound, proved in Lemma 3.5 of \cite{S}:
\begin{equation}\label{bound}
\Delta\rho\geq -\dfrac{\Theta'}{\Theta}\circ\rho
\end{equation}
as distributions on $\Omega$. This means that, if $f$ is any non-negative test-function:
\begin{equation}
(\Delta\rho,f)\geq - \Big(\dfrac{\Theta'}{\Theta}\circ\rho,f\Big)=
-\int_{\Omega}f\Big(\dfrac{\Theta'}{\Theta}\circ\rho\Big).
\end{equation}
We remark that if $\Omega$ is a ball in a space form $M_K$ then the distance function $\rho$ is smooth (except at the center of the ball), the function $\Delta\rho$ is radial (i.e. it depends only on $\rho$) and equality holds in \eqref{bound}.

\smallskip

For $r\geq 0$, we denote by 
$
\{\rho>r\}
$
(respectively, $\{\rho=r\}$) the set of points of $\Omega$ at distance greater (respectively, equal) to $r$ from $\bd\Omega$.
Let $\phi$ be a smooth function on $\Omega$ and let  $F:[0,\infty)\to \reals$ be defined by:
$$
F(r)=\int_{\{\rho>r\}}\phi.
$$
Then $F$ is Lipschitz regular and vanishes for $r\geq R$, where $R$ is the inner radius of $\Omega$. By the formula of co-area, since $\abs{\nabla\rho}=1$ almost everywhere, we have:
$$
F(r)=\int_r^{\infty}\int_{\{\rho=s\}}\phi\, dH_{n-1}\, ds,
$$
where $dH_{n-1}$ is the $(n-1)$-th dimensional Hausdorff measure of the level set $\{\rho=s\}$. Then:
\begin{equation}\label{fprime}
F'(r)=-\int_{\{\rho=r\}}\phi\, dH_{n-1},
\end{equation}
a.e. on $[0,\infty)$. Note that, as $\rho$ is smooth near the boundary, the function $F$ is smooth near $r=0$: precisely on the interval $\big[0,{\rm inj}(\bd\Omega)\big)$, where ${\rm inj}(\bd\Omega)$ is the injectivity radius of the normal exponential map. On that interval one has 
$F'(r)=-\int_{\{\rho=r\}}\phi$ because $dH_{n-1}$ coincides with the Riemannian measure of the smooth hypersurface  $\{\rho=r\}$.

\smallskip

Here is the main result we need. 

\begin{prop}\label{calc} {\rm (\cite{S}, formula (7) p. 517)} 
Let $\phi$ be a smooth function on $\Omega$, and $F(r)=\int_{\{\rho>r\}}\phi$. Then, as distributions on $(0,\infty)$:
\begin{equation}\label{mvl}
F''(r)=-\int_{\{\rho >r\}}\Delta\phi+\rho_{\star}(\phi\Delta\rho),
\end{equation}
where $\rho_{\star}(\phi\Delta\rho)$ denotes the distribution on $(0,\infty)$ given by the push-forward of 
$\phi\Delta\rho$ by $\rho$.
\end{prop}
This means that, if $\psi$ is a test-function on $(0,\infty)$, then, by definition:
$$
(\rho_{\star}(\phi\Delta\rho),\psi)\doteq (\phi\Delta\rho,\psi\circ\rho)=\Big(\Delta\rho,\phi(\psi\circ\rho)\Big),
$$
where on the right hand side $(\cdot,\cdot)$ denotes the duality in $\Omega$. We remark that $\psi\circ\rho$ is only Lipschitz regular; but (see the Appendix below)  $\Delta\rho$ is a zero-order distribution (that is, a Radon measure): then,  the right-hand side is well-defined because it is just the integral of the continuous function $\phi(\psi\circ\rho)$ with respect to the measure $\Delta\rho$. 


\subsection{Proof of Theorem \ref{boundaryintegral}}\label{proofbie}


Fix a non-negative subharmonic function $h$ on $\Omega$, and consider the function
$$
F(r)=\int_{\{\rho>r\}}h.
$$
Let $\psi$ be a non-negative test-function on $(0,\infty)$. Then, by \eqref{bound}, the co-area formula and \eqref{fprime} we see that:
$$
\begin{aligned}
(\rho_{\star}(h\Delta\rho),\psi)&=\Big(\Delta\rho,h(\psi\circ\rho)\Big)\\
&\geq -\Big(\dfrac{\Theta'}{\Theta}\circ\rho,h(\psi\circ\rho)\Big)\\
&= -\int_{\Omega}h(\psi\circ\rho)
\Big(\dfrac{\Theta'}{\Theta}\circ\rho\Big)\\
&=-\int_0^{\infty}\psi(r)\dfrac{\Theta'}{\Theta}(r)\int_{\{\rho=r\}}hdH_{n-1}\,dr\\
&=\int_0^{\infty}\psi(r)\dfrac{\Theta'}{\Theta}(r)F'(r)\,dr\\
&=\Big(\dfrac{\Theta'}{\Theta}F',\psi\Big)
\end{aligned}
$$
that is
\begin{equation}\label{zero}
\rho_{\star}(h\Delta\rho)\geq \dfrac{\Theta'}{\Theta}F'
\end{equation}
as distributions on the half-line. 
 We apply Proposition \ref{calc} to $\phi=h$; as $\Delta h\leq 0$ we conclude by \eqref{zero} that
\begin{equation}\label{final}
F''-\dfrac{\Theta'}{\Theta}F'\geq 0
\end{equation}
as distributions on $(0,\infty)$.

\medskip

Now set  $a=\int_{\Omega}h$ and $b=\int_{\bd\Omega}h$. As $\{\rho>r\}$ is empty for $r\geq R$, we have:
$$
F(0)=a,\quad F'(0)=-b, \quad F(r)=0 \quad\text{if \, $r\geq R$}.
$$
We first assume that $\Theta(R)>0$; we know that $\Theta$ is positive on $[0,R]$ by Proposition \ref{kasue} in the Appendix and then the function
$$
G=\dfrac{F'}{\Theta}
$$
is integrable. By \eqref{final}, the distribution $G'$ is non-negative: it is a simple fact that then $G$ is non-decreasing on a set of full measure in $[0,R]$. As $G$ is regular near $r=0$ and $\Theta(0)=1$ we must have $G(r)\geq G(0)=-b$ a.e. on $[0,R]$, which implies that
\begin{equation}\label{ae}
F'(r)+b\Theta(r)\geq 0
\end{equation}
almost everywhere. We integrate \eqref{ae} on $(0,R)$ to obtain $b\int_0^R\Theta(r)\,dr\geq a$. Finally:
\begin{equation}\label{second}
\dfrac{\int_{\bd\Omega}h}{\int_{\Omega}h}=\dfrac ba\geq\dfrac{1}{\int_0^R\Theta(r)\,dr},
\end{equation} 
which gives the assertion.

\smallskip

Now assume that $\Theta(R)=0$: then, again by Proposition \ref{kasue}, $\Omega$ is a ball of radius $R$ in $M_K$. In that case, the injectivity radius is just the radius of the ball,  the function $F$ is smooth on $[0,R)$ and  the proof carries over as well.
\smallskip

Finally,  assume that $h$ is harmonic and  $\Omega$ is a ball in $M_K$: we have equality in \eqref{bound} and then we have equality at every step of the proof, so that
$$
\dfrac{\int_{\bd\Omega}h}{\int_{\Omega}h}=\dfrac{1}{\int_0^R\Theta(r)dr}.
$$
Note that applying the above to $h=1$ we obtain:
\begin{equation}\label{density}
\dfrac{1}{\int_0^R\Theta(r)dr}
=\isoper.
\end{equation}
In fact, $\Theta$ is the density of the Riemannian measure in normal coordinates around $\bd\Omega$ (see the Appendix below).


\subsection{Proof of Theorem \ref{cheng}}\label{chengproof}


Assume that $M$ is a complete Riemannian manifold with Ricci curvature bounded below by $(n-1)K$, and let $M_K$ be the space form of constant curvature $K$. For any $x_0\in M$  the well-known Bishop-Gromov inequality (see for example \cite{Pe}, Lemma 1.6) states that the function:
$$
V(r)=\dfrac{\abs{B(x_0,r)}}{\abs {B_K(r)}}
$$
is non-increasing on $(0,\infty)$, and tends to $1$ as $r\to 0$ (here $B_K(r)$ is any ball of radius $r$ in $M_K$).  If $r<{\rm inj}(x_0)$ we have $\dfrac d{dr}\abs{B(x_0,r)}=\abs{\bd B(x_0,r)}, \dfrac d{dr}\abs{B_K(r)}=\abs{\bd B_K(r)}$ and as $V'(r)\leq 0$  we obtain:
$$
\dfrac{\abs{\bd B(x_0,r)}}{\abs{B(x_0,r)}}\leq \dfrac{\abs{\bd B_K(r)}}{\abs{B_K(r)}}.
$$
Now, from inequality \eqref{isoper} (which is an equality for geodesic balls in a space form)  we have:
$$
q_1(B(x_0,r))\leq \dfrac{\abs{\bd B(x_0,r)}}{\abs{B(x_0,r)}}\leq \dfrac{\abs{\bd B_K(r)}}{\abs{B_K(r)}}=q_1(B_K(r))
$$
which proves the assertion.


\section{Explicit estimates}\label{ee}


We consider the cases where $K=\{0,-1,1\}$.

\smallskip

First, assume that $\Omega$ is a domain with non-negative Ricci curvature, that is,  $K=0$. One has:
$$
\Theta(r)=(1-Hr)^{n-1}.
$$
The integral $\int_0^R\Theta$ can be explicitly computed, and from Theorem \ref{main} one gets easily:

\begin{thm} \label{eeone} Let $\Omega^{n}$ be a domain with non-negative Ricci curvature and mean curvature bounded below by $H\in\reals$. Let $R$ be the inner radius of $\Omega$. 

\parte a If $H>0$ then
$$
q_1(\Omega)\geq \dfrac{nH}{1-(1-RH)^{n}}.
$$ 

\smallskip

\parte b If $H\geq 0$ then $q_1(\Omega)\geq 1/R$. 

\smallskip

\parte c If $H=-\abs{H}<0$, then 
$$
q_1(\Omega)\geq \dfrac{n\abs{H}}{(1+R\abs{H})^{n}-1}.
$$
\end{thm}

We already observed that a) is sharp. Observe that by Proposition \ref{kasue} one has $1-RH\geq 0$, with equality only for the ball: hence $q_1(\Omega)\geq nH$.

\smallskip

We remark that b) is also sharp. In fact, let $N$ be a closed manifold and $\Omega=N\times (0,2R)$ with the product metric ($\Omega$ is often called a {\it flat cylinder}). Note that the inner radius of $\Omega$ is exactly $R$.  Now $\abs{\bd\Omega}=2\abs{N}$
and $\abs{\Omega}=2R\abs{N}$, so that $\abs{\bd\Omega}/\abs{\Omega}=1/R$. By the upper bound \eqref{isoper} we have $q_1(\Omega)\leq 1/R$ and then, by b), we have equality.

\medskip

We now take $K=-1$, so that the estimates below apply in particular to hyperbolic domains. 

\begin{thm}\label{eetwo} Let $\Omega^n$ be a domain with Ricci curvature bounded below by $-(n-1)$ and mean curvature bounded below by $H\in\reals$. 

\parte a If $H\geq 1$ then $q_1(\Omega)\geq\dfrac{n-1}{1-e^{-(n-1)R}}$. In particular:
$$
q_1(\Omega)> n-1.
$$
\parte b If $H\geq 0$ then $q_1(\Omega)\geq \dfrac{n-1}{e^{(n-1)R}-1}$.

\smallskip

\parte c If $H=-\abs{H}\leq 0$ then $q_1(\Omega)\geq \dfrac{n-1}{(1+\abs{H})^{n-1}(e^{(n-1)R}-1)}.
$
\end{thm} 

For the proof, we  observe that  $K=-1$ and so $\Theta(r)=(\cosh r-H\sinh r)^{n-1}$. To prove a) we use the fact that, if $H\geq 1$, then $\cosh r-H\sinh r\leq e^{-r}$, and then we carry out integration. If $H\geq 0$ we use the inequality $\cosh r \leq e^r$, and if $H=-\abs{H}\leq 0$ we use the inequality $\cosh r-H\sinh r\leq (1+\abs{H})e^r$.

\smallskip

Finally, if the Ricci curvature is bounded below by $n-1$, then 
$
\Theta(r)=(\cos r-H\sin r)^{n-1},
$
and one can get explicit estimates as well. We omit further details.


\subsection{Proof of Corollary \ref{maincor}}\label{proofcor}


Let $\Omega$ be a domain with curvature bounds $(K,H)$ and let $\Theta$ be as in \eqref{theta}.  By Theorem \ref{main} we have:
\begin{equation}\label{compone}
q_1(\Omega)\geq \dfrac{1}{\int_0^R\Theta}.
\end{equation}
Under the given assumptions on $K$ and $H$, there is a unique ball $\bar\Omega$ in $M_K$ with boundary of (constant) mean curvature $H$. By its  definition $\bar\Omega$ has also curvature bounds $(K,H)$  and by the second part of Theorem \ref{main} we  have
\begin{equation}\label{comptwo}
q_1(\bar\Omega)=\dfrac{\abs{\bd\bar\Omega}}{\abs{\bar\Omega}}=
\dfrac{1}{\int_0^{\bar R}\Theta},
\end{equation}
where $\bar R$ is the radius of $\bar\Omega$. Now, by Proposition \ref{kasue},
$\bar R$ is also the first zero of $\Theta$, and one has $\bar R\geq R$ with equality if and only if $\Omega$ is isometric to $\bar\Omega$. Comparison of \eqref{compone} and \eqref{comptwo} leads to the inequality  $q_1(\Omega)\geq q_1(\bar\Omega)$ with equality iff $\Omega$ is isometric to $\bar\Omega$.


\section{Appendix}

We outline the main arguments and definitions in \cite{S} leading to Proposition \ref{calc}.
All these facts are proved in Sections 2 and 3.2 of \cite{S}.

\smallskip

Let $N_x$ be the inner unit normal at $x\in\bd\Omega$, and consider the geodesic normal to $\bd\Omega$ and starting at $x$:
$
\gamma_x(t)={\rm exp}_x(tN_x),
$
where $t$ ranges in a suitable interval.
The {\it cut-radius} of $x\in\bd\Omega$ is the positive number $c(x)$ defined in the following way:

\smallskip

$-$ The geodesic $\gamma_x(t)$ minimizes the distance to $\bd\Omega$ if and only if $t\in[0,c(x)]$. 

\smallskip

The map $c:\bd\Omega\to [0,\infty)$ is continuous; moreover, since $\bd\Omega$ is smooth, $c$ is positive (and $\inf_{\bd\Omega}c$ is called the {\it injectivity radius} of the normal exponential map).
The {\it cut-locus} ${\rm Cut}(\bd\Omega)$ is the closed subset of $\Omega$ defined  by:
$$
{\rm Cut}(\bd\Omega)=\{{\rm exp}_x(c(x)N_x): x\in\bd\Omega\}
$$
It is known that the cut-locus has  measure zero in $\Omega$; denote by
$$
\Omega_{\rm reg}=\Omega\setminus{\rm Cut}(\bd\Omega)
$$
 the set of {\it regular points} of $\rho$. Then, $\rho$ is $C^{\infty}-$smooth on $\Omega_{\rm reg}$ and there one has $\abs{\nabla\rho}=1$. Consider the set
$$
U=\{(r,x)\in [0,\infty)\times\bd\Omega: 0\leq r<c(x)\}.
$$
The pair $(r,x)$ gives rise to the {\it normal coordinates} of a regular point of $\Omega$.
The map $\Phi(r,x)={\rm exp}_x(rN_x)$ is a diffeomorphism
$
\Phi: U\to \Omega_{\rm reg},
$
and if we pull back the Riemannian volume form ${\rm dvol}_g$ of $\Omega$ by $\Phi$,  we can write
$$
\starred\Phi({\rm dvol}_g)(r,x)=\theta(r,x)dr\,dx,
$$
where $dx$ denotes, for short,  the induced volume form of $\bd\Omega$. The function $\theta$ is then the density of the Riemannian volume form in normal coordinates. Obviously $\theta$ is positive on $U$, and $\theta(0,x)=1$ for all $x\in\bd\Omega$. So, for all integrable functions $f$ on $\Omega$:
\begin{equation}\label{nc}
\int_{\Omega}f=\int_{\Omega_{\rm reg}}f=\int_{\bd\Omega}\int_0^{c(x)}\theta(r,x) f(r,x)\,dr\,dx
\end{equation}
where we identify a regular point of $\Omega$ with its pair of normal coordinates. 

The map $\theta$ extends by continuity on $\bar U$, and we will define
$$
\theta(c(x),x)\doteq\lim_{r\to c(x)_-}\theta(r,x).
$$
We let
$$
\Delta_{\rm reg}\rho=\Delta (\rho|_{\Omega_{\rm reg}}),
$$
be the {\it regular part} of the Laplacian of the distance function. It is an $L^1$- function on $\Omega_{\rm reg}$. In normal coordinates, one has the formula (see \cite{Ga} p. 40):
\begin{equation}\label{reg}
\Delta_{\rm reg}\rho(r,x)=-\dfrac{1}{\theta}\derive{\theta}{r}(r,x).
\end{equation}
Geometrically, $\Delta_{\rm reg}\rho(r,x)$ is equal to $(n-1)$-times the mean curvature of the level set $\{\rho=r\}$ at the regular point $(r,x)$.

\medskip

In what follows, we assume that $\Omega$ has curvature bounds $(K,H)$; that is, its Ricci curvature is bounded below by $(n-1)K$ and its boundary $\bd\Omega$ has mean curvature bounded below by $H$.
Recall that 
$
\Theta(r)=(s'_K(r)-Hs_K(r))^{n-1},
$
where $s_K(r)$ has been defined in \eqref{sk}. Then, the classical volume estimates of Heintze and Karcher (for an independent derivation using a Laplacian comparison argument see p. 41 of \cite{Ga})
imply that, at all regular points $(r,x)\in U$ one has
\begin{equation}\label{regular}
\Delta_{\rm reg}\rho(r,x)=-\dfrac{1}{\theta}\derive{\theta}{r}(r,x)\geq -\dfrac{\Theta'(r)}{\Theta(r)}
\end{equation}
This implies in particular that, on $U$, one has 
$\theta(r,x)\leq\Theta(r)$. As a consequence of this, and a detailed analysis of the focal points of $\bd\Omega$, one can prove (see Theorem A of \cite{Ka}):

\begin{prop}\label{kasue} 
Assume that $\Omega$ has curvature bounds $(K,H)$, and let $\bar R\in (0,\infty]$ be the first positive zero of the function $r\to s'_K(r)-Hs_K(r)$. Then
$$
\bar R\geq R,
$$
where $R$ is the inner radius of $\Omega$. Moreover, equality holds if and only if $\Omega$ is a ball of radius $R$ in the space form $M_K$. 
Then, the function $\Theta$ is smooth and positive on the interval $[0,R)$ and $\Theta (R)=0$ only when $\Omega$ is a geodesic ball in $M_K$.
\end{prop} 

As a distribution on $\Omega$, the Laplacian $\Delta\rho$ splits as follows:
\begin{equation}\label{splitting}
\Delta\rho=\Delta_{\rm reg}\rho+\Delta_{\rm cut}\rho,
\end{equation}
where the regular part $\Delta_{\rm reg}\rho$ is  $L^1$ and satisfies \eqref{regular} and where
$\Delta_{\rm cut}\rho$ is the distribution supported on the cut-locus and defined by:
$$
(\Delta_{\rm cut}\rho,f)=\int_{\bd\Omega}\theta(c(x),x)f(c(x),x)\,dx.
$$
Note that  $\Delta_{\rm cut}\rho$ is a non-negative distribution: as such it can be identified with a non-negative Radon measure on $\Omega$.  Hence $\Delta\rho$ is a (signed) Radon measure on $\Omega$, and can be tested on any continuous function. 
To verify the splitting in \eqref{splitting} (proved in  Lemma 2.1 and Section 3.2 of \cite{S}), we test $\Delta\rho$ on a smooth function, use normal coordinates, and integrate by parts. 
From the above splitting, the positivity  of $\Delta_{\rm cut}\rho$ and inequality \eqref{regular} one then obtains  the bound $
\Delta\rho\geq -\frac{\Theta'}{\Theta}\circ\rho$ in \eqref{bound}. 

Finally, the proof of Proposition \ref{calc} is done using integrating by parts on the level domains $\{\rho>r\}$. We remark that the arguments in \cite{S} extend to the distance function to any submanifold of a Riemannian manifold. 

\begin{rem}\label{tball} 
\rm When $\Omega$ is a ball of radius $R$ in $M_K$, one sees that, by the symmetries of $M_K$, the density $\theta(r,x)$ is independent on $x\in\bd\Omega$ and in fact one has $\theta(r,x)=\Theta(r)$ for all $r,x$; moreover $c(x)= R$ for all $x\in\bd\Omega$, and $\Theta( R)=0$. The cut-locus reduces to a point; this is the unique focal point of $\bd\Omega$ and  coincides with the center of the ball: hence $\Delta_{\rm cut}\rho=0$. We have equality in \eqref{regular} hence also in \eqref{bound} and from \eqref{nc} one verifies that
\begin{equation}\label{ball}
\dfrac{1}{\int_0^R\Theta(r)dr}=\dfrac{\abs{\bd\Omega}}{\abs{\Omega}}.
\end{equation}
\end{rem}


\vspace{0.8cm}     
Authors addresses:     
\nopagebreak     
\vspace{5mm}\\     
\parskip0ex     
\vtop{\hsize=6cm\noindent\obeylines}     
\vtop{     
\hsize=8cm\noindent     
\obeylines     
Simon Raulot
Laboratoire de Math\'ematiques R. Salem
UMR $6085$ CNRS-Universit\'e de Rouen
Avenue de l'Universit\'e, BP.$12$
Technop\^ole du Madrillet
$76801$ Saint-\'Etienne-du-Rouvray, France}     
     
\vspace{0.5cm}     
     
E-Mail:     
{\tt simon.raulot@univ-rouen.fr }  

\vtop{\hsize=6cm\noindent\obeylines}     
\vtop{     
\hsize=9cm\noindent     
\obeylines     
Alessandro Savo
Dipartimento SBAI, Sezione di Matematica 
Sapienza Universit\`a di Roma
Via Antonio Scarpa 16
 00161 Roma, Italy         
}     
     
\vspace{0.5cm}     
     
E-Mail:     
{\tt alessandro.savo@sbai.uniroma1.it  } 



\end{document}